\documentclass[a4paper,12pt]{article}

 \usepackage{latexsym,amssymb,amsmath}
 
 \newtheorem{thm}{Theorem}
 \newtheorem{lem}{Lemma}

\begin{document}

\title{Shifted Appell sequences in Clifford analysis}

\author{Dixan Pe\~na Pe\~na\\\normalsize{Department of Mathematics, University of Aveiro}\\\normalsize{3810-193 Aveiro, Portugal}\\\normalsize{e-mail: dixanpena@ua.pt; dixanpena@gmail.com}}

\date{} 

\maketitle

\begin{abstract}
\noindent This paper is a continuation of [D. Pe\~{n}a Pe\~{n}a, On a sequence of monogenic polynomials satisfying the Appell condition whose first term is a non-constant function, arXiv:1102.1833], in which we prove that for every monogenic polynomial $\mathbf{P}_k(x)$ of degree $k$ in $\mathbb R^{m+1}$ there exists a sequence of monogenic polynomials $\{M_n(x)\}_{n\ge0}$ satisfying the Appell condition such that $M_0(x)=\mathbf{P}_k(x)$.\vspace{0.2cm}\vspace{0.1cm}
\end{abstract}

\section{Preliminaries}

The real Clifford algebra $\mathbb R_{0,m}$ (see \cite{Cl}) is the free algebra generated by  the standard basis $\{e_1,\ldots,e_m\}$ of the Euclidean space $\mathbb R^m$, subject to the multiplication relations 
\[e_je_k+e_ke_j=-2\delta_{jk},\quad j,k=1,\dots,m,\]
where $\delta_{jk}$ denotes the Kronecker delta. The dimension of the real Clifford algebra $\mathbb R_{0,m}$ is $2^m$, as is the case for the Grassmann algebra generated by $\{e_1,\ldots,e_m\}$, but the difference is that now $e_j^2=-1$ instead of $e_j^2=0$, creating a structure with similarities to the complex numbers. 

A general element $a$ of $\mathbb R_{0,m}$ may be written as $a=\sum_Aa_Ae_A$, $a_A\in\mathbb R$, in terms of the basic elements $e_A=e_{j_1}\dots e_{j_k}$, defined for every subset $A=\{j_1,\dots,j_k\}$ of $\{1,\dots,m\}$ with $j_1<\dots<j_k$. For the empty set, one puts $e_{\emptyset}=1$, the latter being the identity element. Conjugation in $\mathbb R_{0,m}$ is given by $\overline a=\sum_Aa_A\overline e_A$, with $\overline e_A=\overline e_{j_k}\dots\overline e_{j_1}$, $\overline e_j=-e_j$, $j=1,\dots,m$.

One natural way to generalize the holomorphic functions to higher dimensions is by considering the null solutions of the fundamental first order differential operator $\partial_x$ in $\mathbb R^{m+1}$ given by 
\[\partial_x=\partial_{x_0}+\partial_{\underline x}=\partial_{x_0}+\sum_{j=1}^me_j\partial_{x_j},\]
called the generalized Cauchy-Riemann operator, and where $\partial_{\underline x}$ is the Dirac operator in $\mathbb R^m$. That is,  an $\mathbb R_{0,m}$-valued function $f$ defined and continuously differentiable in an open set $\Omega$ of $\mathbb R^{m+1}$, is said to be (left) monogenic in $\Omega$ if and only if $\partial_xf=0$ in $\Omega$. In a similar way one also defines monogenicity with respect to the Dirac operator $\partial_{\underline x}$ in $\mathbb R^m$. Monogenic functions are a central object of study in Clifford analysis (see e.g. \cite{BDS,DSS,GuSp}).   

An remarkable feature of the generalized Cauchy-Riemann operator $\partial_x$ is that it gives a factorization of the Laplacian, i.e.
\[\Delta_x=\sum_{j=0}^m\partial_{x_j}^2=\partial_x\overline\partial_x=\overline\partial_x\partial_x,\]
and therefore every monogenic function is also harmonic. Observe that the operator $\overline\partial_x=\partial_{x_0}-\partial_{\underline x}$ may be seen as the higher dimensional version of the well-known operator $2\partial_z=\partial_{x_0}-i\partial_{x_1}$. Furthermore, according to \cite{GuM,M}, the hypercomplex derivative of a monogenic function $f$ is defined as 
\[\frac{1}{2}\,\overline\partial_xf.\] 
As a monogenic function $f$ clearly satisfies $\partial_{x_0}f=-\partial_{\underline x}f$, it easily follows that
\[\frac{1}{2}\,\overline\partial_xf=\partial_{x_0}f=-\partial_{\underline x}f.\]
An important class of polynomial sequences is the class of Appell sequences which is defined as follows (see \cite{A}). A polynomial sequence $\{p_n(t)\}_{n\ge0}$, i.e. the index of each polynomial equals its degree, is said to be an Appell sequence if it satisfies
\[p_n^\prime(t)=np_{n-1}(t),\quad n\ge1.\]
Probably, the simplest example of an Appell sequence is the sequence $\{t^n\}_{n\ge0}$, other examples being the Bernoulli, the Euler and the Hermite polynomials. 

Appell sequences have been recently introduced to the Clifford analysis setting (see e.g. \cite{BG,BGLS,CM,BieS,FCM,FM,NGu,Lav,MF}). Namely, a sequence $\{P_n(x)\}_{n\ge0}$ of $\mathbb R_{0,m}$-valued polynomials forms an Appell sequence if the following conditions are satisfied:
\begin{itemize}
\item [{\rm (i)}] $\{P_n(x)\}_{n\ge0}$ is a polynomial sequence;
\item [{\rm (ii)}] each $P_n(x)$ is monogenic in $\mathbb R^{m+1}$, i.e. $\partial_xP_n(x)=0$ for all $x\in\mathbb R^{m+1}$;
\item [{\rm (iii)}] $\frac{1}{2}\,\overline\partial_xP_n(x)=nP_{n-1}(x)$, $n\ge1$.
\end{itemize}
Note that the requirement of $\{P_n(x)\}_{n\ge0}$ being a polynomial sequence implies that the first term $P_0(x)$ must be a constant. It is natural to ask whether one can construct sequences of monogenic polynomials satisfying the Appell condition (iii) but in which the first term is a monogenic polynomial in $\mathbb R^{m+1}$ and not necessarily a constant. More precisely, we are interested in sequences $\{M_n(x)\}_{n\ge0}$ of $\mathbb R_{0,m}$-valued polynomials which are monogenic in $\mathbb R^{m+1}$ fulfilling
\begin{equation}\label{Acond}
\frac{1}{2}\,\overline\partial_xM_n(x)=nM_{n-1}(x),\quad n\ge1,
\end{equation}
where $M_0(x)$ is an arbitrary monogenic polynomial in $\mathbb R^{m+1}$. These sequences will be called \emph{shifted Appell sequences of monogenic polynomials}. 

This paper is a continuation of \cite{DPP}, where an example of these sequences was constructed for the case $M_0(x)=\mathbf{P}_k(\underline x)$ being an arbitrary $\mathbb R_{0,m}$-valued homogeneous monogenic polynomial of degree $k$ in $\mathbb R^m$. 

It is clear that the class of shifted Appell sequences of monogenic polynomials is a right $\mathbb R_{0,m}$-module under the usual addition of sequences and multiplication by Clifford numbers. Suppose now that $\mathbf{P}_\kappa(x)$ is an $\mathbb R_{0,m}$-valued polynomial of degree $\kappa$ which moreover is monogenic in $\mathbb R^{m+1}$. It is easy to check that $\mathbf{P}_\kappa(x)$ may be written as
\[\mathbf{P}_\kappa(x)=\sum_{k=0}^\kappa\mathbf{P}_k(x),\]
where $\mathbf{P}_k(x)$ denotes a homogeneous monogenic polynomial of degree $k$ in $\mathbb R^{m+1}$. Thus, on account of the previous remarks, we only need to prove:

\begin{thm}\label{mainresult}
Let $\mathbf{P}_k(x)$ be an $\mathbb R_{0,m}$-valued homogeneous polynomial of degree $k$ which is monogenic in $\mathbb R^{m+1}$. Then there exists a shifted Appell sequence of monogenic polynomials $\{M_n(x)\}_{n\ge0}$ such that $M_0(x)=\mathbf{P}_k(x)$.  
\end{thm}

\section{Some fundamental results}

Let $\mathsf{P}(k)$ ($k\in\mathbb N_0$) be the set of all $\mathbb R_{0,m}$-valued homogeneous polynomials of degree $k$ in $\mathbb R^m$. This set contains the important subspace $\mathsf{M}^+(k)$ consisting of all polynomials in $\mathsf{P}(k)$ which are monogenic. That is, $\mathbf{P}_k(\underline x)\in\mathsf{M}^+(k)$ if it is an $\mathbb R_{0,m}$-valued polynomial of degree $k$ and
\[\mathbf{P}_k(t\underline x)=t^k\mathbf{P}_k(\underline x),\quad\partial_{\underline x}\mathbf{P}_k(\underline x)=0,\quad \underline x\in\mathbb R^m,\;t\in\mathbb R.\]
For a differentiable $\mathbb R$-valued function $\phi$ and a differentiable $\mathbb R_{0,m}$-valued function $g$, we have
\begin{equation}\label{Lr1}
\partial_{\underline x}(\phi g)=\partial_{\underline x}(\phi)g + \phi(\partial_{\underline x}g).
\end{equation}
Moreover, for a differentiable vector-valued function $\underline f=\sum_{j=1}^mf_je_j$, we also have
\begin{equation}\label{Lr2}
\partial_{\underline x}(\underline fg)=(\partial_{\underline x}\underline f)g-\underline f(\partial_{\underline x}g)-2\sum_{j=1}^mf_j(\partial_{x_j}g).
\end{equation}
Let  
\[\beta_k(n)=\left\{\begin{array}{ll}n,&\text{if}\;n\;\text{even}\\2k+m+n-1,&\text{if}\;n\;\text{odd}\end{array}\right.\]
for $n\ge1$. Using the Leibniz rules (\ref{Lr1})-(\ref{Lr2}) as well as Euler's theorem for homogeneous functions, we can deduce the useful identity:
\begin{equation}\label{ident1}
\partial_{\underline x}\big(\underline x^n\mathbf{P}_k(\underline x)\big)=-\beta_k(n)\underline x^{n-1}\mathbf{P}_k(\underline x),\quad\mathbf{P}_k(\underline x)\in\mathsf{M}^+(k),\quad n\ge1.
\end{equation}
Let us recall two basic results of Clifford analysis: the Cauchy-Kovalevskaya extension technique (see e.g. \cite{BDS,DSS}) and the Almansi-Fischer decomposition (see e.g. \cite{DSS,MR}). 

\begin{thm}\label{CK}
Every $\mathbb R_{0,m}$-valued function $g(\underline x)$ analytic in $\mathbb R^m$ has a unique monogenic extension $\mathsf{CK}[g]$ to $\mathbb R^{m+1}$, which is given by
\begin{equation*}\label{CKf}
\mathsf{CK}[g(\underline x)](x)=\sum_{j=0}^\infty\frac{(-x_0)^j}{j!}\,\partial_{\underline x}^jg(\underline x).
\end{equation*}
\end{thm}
\textbf{Remark:} Observe that a monogenic function $f(x)$ can be reconstructed by knowing its restriction to $\mathbb R^m$ using previous formula, i.e. \[f(x)=\mathsf{CK}[f(x)\vert_{x_0=0}](x).\] 
It is also worth noting that
\begin{equation}\label{DCK}
\frac{1}{2}\,\overline\partial_x\mathsf{CK}[g(\underline x)](x)=-\partial_{\underline x}\mathsf{CK}[g(\underline x)](x)=-\mathsf{CK}[\partial_{\underline x}g(\underline x)](x).
\end{equation}

\begin{thm}\label{Fischer}
Let $k\in\mathbb N$. Then
\[\mathsf{P}(k)=\bigoplus_{\nu=0}^k\underline x^\nu\mathsf{M}^+(k-\nu).\]
\end{thm}
Theorems \ref{CK} and \ref{Fischer} together with equality (\ref{ident1}) will be essential for proving our main result.

\section{Proof of the main result}

We shall first introduce a collection $\big\{\{\mathsf{M}_n^{k,\nu}(x)\}_{n\ge0},\;0\le\nu\le k\big\}$ of shifted Appell sequences of homogeneous monogenic polynomials whose first terms are  
\[\mathsf{M}_0^{k,\nu}(x)=\mathsf{CK}[\underline x^\nu\mathbf{P}_{k-\nu}(\underline x)](x),\quad\mathbf{P}_{k-\nu}(\underline x)\in\mathsf{M}^+(k-\nu),\quad0\le\nu\le k.\]  
From Theorem \ref{CK} we can deduce that $\mathsf{M}_0^{k,\nu}(x)$ is a homogeneous monogenic polynomials of degree $k$ in $\mathbb R^{m+1}$ and is of the form  
\[\mathsf{M}_0^{k,\nu}(x)=\left(\sum_{j=0}^\nu\frac{\mu_j^{k,\nu}}{j!}x_0^j\underline x^{\nu-j}\right)\mathbf{P}_{k-\nu}(\underline x),\]
with $\mu_j^{k,\nu}=\displaystyle{\prod_{s=\nu-j+1}^\nu\beta_{k-\nu}(s)}$ for $1\le j\le\nu$ and $\mu_j^{k,\nu}=1$ for $j=0$.

\begin{lem}\label{lemainter}
Assume that $\mathbf{P}_{k-\nu}(\underline x)\in\mathsf{M}^+(k-\nu)$ where $k,\nu\in\mathbb N_0$, $\nu\le k$ and put 
\[\lambda_n^{k,\nu}=\frac{n!}{\displaystyle{\prod_{s=1}^{n}\beta_{k-\nu}(\nu+s)}},\]
for $n\ge1$ and $\lambda_n^{k,\nu}=1$ for $n=0$. The sequence $\{\mathsf{M}_n^{k,\nu}(x)\}_{n\ge0}$ defined by 
\[\mathsf{M}_n^{k,\nu}(x)=\lambda_n^{k,\nu}\mathsf{CK}[\underline x^{\nu+n}\mathbf{P}_{k-\nu}(\underline x)](x),\quad n\ge0,\]
is a shifted Appell sequence of homogeneous monogenic polynomials.
\end{lem}
\textit{Proof.} By Theorem \ref{CK}, it follows that $\mathsf{M}_n^{k,\nu}(x)$ is a homogeneous monogenic polynomials of degree $k+n$ in $\mathbb R^{m+1}$ and is of the form
\[\mathsf{M}_n^{k,\nu}(x)=n!\left(\sum_{j=0}^{n-1}\frac{\lambda_{n-j}^{k,\nu}}{j!(n-j)!}x_0^j\underline x^{\nu+n-j}+\sum_{j=n}^{\nu+n}\frac{\mu_{j-n}^{k,\nu}}{j!}x_0^j\underline x^{\nu+n-j}\right)\mathbf{P}_{k-\nu}(\underline x),\,n\ge1.\]
It only remains to show that $\{\mathsf{M}_n^{k,\nu}(x)\}_{n\ge0}$ satisfies the Appell condition (\ref{Acond}). Indeed, using (\ref{DCK}) and identity (\ref{ident1}), we see at once that 
\begin{align*}
\frac{1}{2}\,\overline\partial_x\mathsf{M}_n^{k,\nu}(x)&=-\lambda_n^{k,\nu}\mathsf{CK}\big[\partial_{\underline x}\big(\underline x^{\nu+n}\mathbf{P}_{k-\nu}(\underline x)\big)\big](x)\\
&=\beta_{k-\nu}(\nu+n)\lambda_n^{k,\nu}\mathsf{CK}\big[\underline x^{\nu+n-1}\mathbf{P}_{k-\nu}(\underline x)\big](x)\\
&=n\lambda_{n-1}^{k,\nu}\mathsf{CK}[\underline x^{\nu+n-1}\mathbf{P}_{k-\nu}(\underline x)](x)=n\mathsf{M}_{n-1}^{k,\nu}(x).\quad\square
\end{align*}
\textbf{Remark:} It should be noticed that $\{\mathsf{M}_n^{k,0}(x)\}_{n\ge0}$, which is the first sequence in the above collection, corresponds to the sequence constructed in \cite{DPP}. 

\noindent We can now prove the main result of this paper:

\noindent\textit{Proof of Theorem \ref{mainresult}.} Suppose that $\mathbf{P}_k(x)$ is a $\mathbb R_{0,m}$-valued homogeneous polynomial of degree $k$ which is monogenic in $\mathbb R^{m+1}$. From Theorem \ref{CK} we have that 
\[\mathbf{P}_k(x)=\mathsf{CK}\left[\mathbf{P}_k(x)\vert_{x_0=0}\right](x).\]
It is clear that $\mathbf{P}_k(x)\vert_{x_0=0}\in\mathsf{P}(k)$. Consequently, by Theorem \ref{Fischer}, there exists unique $\mathbf{P}_{k-\nu}(\underline x)\in\mathsf{M}^+(k-\nu)$ such that
\[\mathbf{P}_k(x)\vert_{x_0=0}=\sum_{\nu=0}^k\underline x^\nu\mathbf{P}_{k-\nu}(\underline x).\]
From the above it follows that
\[\mathbf{P}_k(x)=\sum_{\nu=0}^k\mathsf{CK}\left[\underline x^\nu\mathbf{P}_{k-\nu}(\underline x)\right](x).\]
Define
\[\{M_n(x)\}_{n\ge0}=\left\{\sum_{\nu=0}^k\mathsf{M}_n^{k,\nu}(x)\right\}_{n\ge0}.\]
Lemma \ref{lemainter} now shows that $\{M_n(x)\}_{n\ge0}$ is a shifted Appell sequence of monogenic polynomials with first term $M_0(x)=\mathbf{P}_k(x).\quad\square$ 

\subsection*{Acknowledgment}

The author was supported by a Post-Doctoral Grant of \emph{Funda\c{c}\~ao para a Ci\^encia e a Tecnologia} (FCT), Portugal.

\end{document}